# Realisation of special Kähler manifolds as parabolic spheres

Oliver Baues *  
Departement Mathematik  
ETH-Zentrum  
Rämistrasse 101  
CH-8092 Zürich

Vicente Cortés †  
Mathematisches Institut  
Universität Bonn  
Beringstraße 1  
D-53115 Bonn

October 21, 1999

**Abstract**

We prove that any simply connected special Kähler manifold admits a canonical immersion as a parabolic affine hypersphere. As an application, we associate a parabolic affine hypersphere to any nondegenerate holomorphic function. Also we show that a classical result of Calabi and Pogorelov on parabolic spheres implies Lu's theorem on complete special Kähler manifolds with a positive definite metric.

## Introduction

The purpose of this note is to relate some area of recent intense interest, namely the special Kähler manifolds, to the more classical subject of affine differential geometry. We show that any simply connected special Kähler manifold admits a canonical immersion as a parabolic affine hypersphere (Theorem 3.1). Also we characterise the parabolic affine hyperspheres which arise in this way (Theorem 3.6). The results show that the special Kähler condition is in fact equivalent to the requirement that a Kähler manifold locally embeds into affine space as a parabolic sphere in a compatible way.

The correspondence between special Kähler manifolds and a subclass of parabolic affine hyperspheres allows to relate results from both areas of research. For example, special Kähler manifolds are known to have particular good local properties. All of them may be locally described by nondegenerate holomorphic functions. Blaschke observed (see [B]) that the local equations of a two dimensional parabolic sphere are integrable with the help of a holomorphic function. Since in dimension two every Riemannian parabolic sphere is obtained locally as an immersion of a special Kähler manifold, our correpondence generalises Blaschke's result to higher dimensions. As another application, we remark that a recent theorem of Lu on complete special Kähler manifolds with a positive definite metric (see [L]) reduces to a special case of a well known result of Calabi and Pogorelov on parabolic spheres.

*e-mail: oliver@math.ethz.ch  
†e-mail: vicente@math.uni-bonn.de

This work was supported by SFB256 (Universität Bonn)

The paper is organized as follows: In section 1 we recall some basic facts and definitions of affine differential geometry. In particular, we introduce the so called *Fundamental Theorem* which is our basic tool. In section 2 we recall the definition of special Kähler manifolds and explain their local structure. In section 3 we prove our main theorems and give the two applications mentioned above.

# 1 Affine hyperspheres

Let us begin by recalling some basic facts about affine immersions. For a comprehensive treatment of the subject, see the monograph [NS]. We consider $\mathbf{R}^{n+1}$ as affine space with standard volume form det. An affine transformation is called *equiaffine* if it preserves det. Let $\varphi : M \to \mathbf{R}^{n+1}$ be an immersion of an $n$-dimensional manifold $M$ into affine $(n+1)$-space. The choice of a transversal vector field along $M$ induces the following data on the tangent bundle of $M$: a torsionfree connection $\nabla$, a symmetric bilinear form $g$, a volume form $\theta = \det(\xi, \ldots)$, and an endomorphism $S$ called the shape tensor.

If $g$ is nondegenerate (a condition which is independent of the choice of transversal field) then there exists a canonical transversal field $\xi$, called the *affine normal*. The affine normal $\xi$ is uniquely determined up to sign by the following conditions: The metric volume form $\nu$ of $g$ is $\nabla$-parallel and coincides with the induced volume form $\theta$. The immersion $\varphi$ with the affine normal $\xi$ is called a *Blaschke immersion*, and the pseudo-Riemannian metric $g$ associated to the affine normal is called the *Blaschke metric* of the nondegenerate hypersurface.

The pair $(\nabla, g)$ associated to the immersion $\varphi$ satisfies certain integrability conditions. The fundamental theorem of affine differential geometry asserts that given a simply connected manifold $M$ with (abstract) data $(\nabla, g)$ satisfying the integrability conditions, there exists an immersion $\varphi : M \to \mathbf{R}^{n+1}$ inducing the pair $(\nabla, g)$, see [NS, p. 21 and p. 75]. Before we state the theorem, let us recall that the equation:

$$Xg(Y, Z) = g(\nabla_X Y, Z) + g(Y, \overline{\nabla}_X Z)$$

defines the $g$-conjugate connection $\overline{\nabla}$ on $M$.

**Theorem 1.1** *Let $M$ be a simply connected manifold with a torsionfree connection $\nabla$ and with a pseudo-Riemannian metric $g$. Then there exists a Blaschke immersion $\varphi : M \to \mathbf{R}^{n+1}$, unique up to equiaffine transformations, with Blaschke metric $g$ and induced connection $\nabla$ if and only if the $g$-conjugate connection $\overline{\nabla}$ is torsionfree, projectively flat, and if the metric volume form $\nu$ is $\nabla$-parallel.*

We note that the shape tensor of a Blaschke immersion can be expressed in terms of $g$ and the Ricci curvature $\mathrm{Ric}^{\overline{\nabla}}$ of the conjugate connection $\overline{\nabla}$ as follows:

$$g(SX, Y) = \mathrm{Ric}^{\overline{\nabla}}(X, Y).$$

A Blaschke immersion is called an *affine hypersphere* if its shape tensor satisfies $S = \lambda \, \mathrm{Id}$, for some $\lambda \in \mathbf{R}$. An affine hypersphere is called *proper* if $\lambda \neq 0$ and it is called *parabolic* (or improper) if $\lambda = 0$. Later we will use the following characterisation of parabolic spheres, see [NS, p. 42].

**Proposition 1.2** *A Blaschke immersion $\varphi : M \to \mathbf{R}^{n+1}$ is a parabolic affine sphere if and only if the induced connection $\nabla$ is flat.*

## 2 Special Kähler manifolds

Let us first recall the definition of an (affine) special Kähler manifold, see e.g. [F]. Let $(M,g)$ be a Kähler manifold with, possibly indefinite, Kähler metric $g$ and with Kähler form $\omega := g(\cdot, J\cdot)$, where $J$ denotes the complex structure on $M$. Assume also that there is given a torsionfree, flat connection $\nabla$ on $M$. The triple $(M, \nabla, g)$ is called a *special Kähler manifold* if $\nabla \omega = 0$ and $d^\nabla J = 0$, where $d^\nabla J(X, Y) = (\nabla_X J)Y - (\nabla_Y J)X$.

We consider now the complex vector space $V = T^*\mathbf{C}^m$ with standard complex symplectic form $\Omega = \sum_{i=1}^m dz^i \wedge dw_i$. Let $\tau : V \to V$ be the standard real structure of $V$ with set of fixed points $V^\tau = T^*\mathbf{R}^m$. Then $\gamma := \sqrt{-1}\,\Omega(\cdot, \tau\cdot)$ defines a Hermitian form of (complex) signature $(m, m)$. A holomorphic immersion $\phi : M \to V$ of a complex manifold $M$ into $V$ is called *nondegenerate* if $\phi^*\gamma$ is nondegenerate. If $\phi$ is nondegenerate $\phi^*\gamma$ defines a, possibly indefinite, Kähler metric $g$ on $M$. If, additionally, $\phi$ is a Lagrangian immersion then it induces a torsionfree flat connection $\nabla$ on $M$, see [C]. It is proven in [ACD] that $(M, \nabla, g)$ is a special Kähler manifold.

It is well known from classical mechanics that any holomorphic Lagrangian immersion $\phi$ may be locally represented as a closed holomorphic 1-form

$$\alpha : U \to T^*\mathbf{C}^m$$

on a domain $U \subset \mathbf{C}^m$. Since we can assume that $U$ is simply connected $\phi$ is locally defined by a holomorphic function $F$ on $U$ with $\alpha = dF$. The holomorphic function $F$ will be called *nondegenerate* if the Lagrangian immersion $dF$ is nondegenerate. Note that in terms of the canonical coordinates $(z^1, \ldots, z^m, w_1, \ldots, w_m)$ the image of $dF$ is defined by the equations $w_i = \partial F/\partial z^i$.

The local structure of a special Kähler manifold is described by the following theorem, see [ACD], which is an analogous result to Theorem 1.1. It shows that any simply connected special Kähler manifold is obtained by the above construction, and hence is determined by a nondegenerate holomorphic function $F$.

**Theorem 2.1** *Let $(M, \nabla, g)$ be a simply connected special Kähler manifold of complex dimension $m$. Then there exists a nondegenerate holomorphic Lagrangian immersion $\phi : M \to V = T^*\mathbf{C}^m$ which induces $\nabla$ and $g$ on $M$. The immersion $\phi$ is uniquely determined up to affine transformations of $V$ preserving the symplectic form $\Omega$ and the real structure $\tau$.*

## 3 Affine immersions of Kähler manifolds

Let $(M, g)$ be a Kähler manifold with Kähler form $\omega$. Given a torsionfree connection $\nabla$ such that $\nabla \omega = 0$, we may ask if there exists a Blaschke immersion with induced data $(\nabla, g)$. Our main result is that any simply connected special Kähler manifold admits such a realisation. Moreover, if $\nabla$ is flat and there exists a Blaschke immersion of $(M, g)$ then $(M, \nabla, g)$ is special Kähler. Thus, we show that the special Kähler condition $d^\nabla J = 0$ is equivalent to the existence of a Blaschke immersion of $(M, g, \nabla)$ under the assumption that $\nabla$ is flat.

**Theorem 3.1** *Let $(M, \nabla, g)$ be a simply connected special Kähler manifold of dimension $n = 2m$. Then there exists a Blaschke immersion $\varphi : M \to \mathbf{R}^{n+1}$ with induced connection $\nabla$ and Blaschke metric $g$. Moreover, $\varphi$ is a parabolic sphere and is uniquely determined up to equiaffine transformations.*

*Proof* We show that the pair $(\nabla, g)$ satisfies the integrability conditions of Theorem 1.1. The metric volume form $\nu$ is $\nabla$-parallel because the Kähler form $\omega$ is $\nabla$-parallel and $\nu$ is a constant multiple of $\omega^m$. We prove now that the $g$-conjugate connection $\overline{\nabla}$ is torsionfree and flat, and hence is also projectively flat. We consider, as in [ACD], the connection $\nabla^J$ defined by $\nabla^J := J\nabla J^{-1}$. Since $\nabla$ is flat, $\nabla^J$ is flat as well. One can easily check that, given a torsionfree flat connection $\nabla$ on a complex manifold $M$, the special condition $d^\nabla J = 0$ is equivalent to the vanishing of the torsion of $\nabla^J$. Now the vanishing of the torsion and curvature of $\overline{\nabla}$ follows from the fact that $\overline{\nabla}$ coincides with the torsionfree and flat connection $\nabla^J$. In fact, the next computation, in which we use only the fact that $\omega$ is parallel and Hermitian, shows that $\nabla^J$ is $g$-conjugate to $\nabla$:

$$\begin{aligned} Xg(Y,Z) &= X\omega(Y, JZ) = \omega(\nabla_X Y, JZ) + \omega(Y, \nabla_X JZ) \\ &= g(\nabla_X Y, Z) + \omega(JY, J\nabla_X JZ) = g(\nabla_X Y, Z) - \omega(JY, \nabla_X^J Z) \\ &= g(\nabla_X Y, Z) + \omega(\nabla_X^J Z, JY) = g(\nabla_X Y, Z) + g(\nabla_X^J Z, Y) \\ &= g(\nabla_X Y, Z) + g(Y, \nabla_X^J Z)\,, \end{aligned}$$

where $X$, $Y$ and $Z$ are vector fields on $M$. Finally, we remark that, since $\nabla$ is flat, $\varphi$ defines a parabolic sphere due to Proposition 1.2. $\square$

The following theorem is a classical result of Calabi and Pogorelov, see [NS, p. 125].

**Theorem 3.2** *If the Blaschke metric $g$ of a parabolic sphere $M$ is (positive) definite and complete, then $M$ is affinely congruent to the paraboloid $x^{n+1} = \sum_{i=1}^n (x^i)^2$ in $\mathbf{R}^{n+1}$. In particular, $g$ is flat.*

As a corollary we obtain the following recent theorem of Lu [L].

**Corollary 3.3** *Let $(M, \nabla, g)$ be a special Kähler manifold with (positive) definite metric $g$. If $g$ is complete then it is flat.*

*Proof* Passing to the universal covering we can assume that $M$ is simply connected, and hence, according to Theorem 3.1, admits a Blaschke immersion, which is a parabolic sphere and which induces $\nabla$ and $g$. Now the corollary is an immediate consequence of Theorem 3.2. $\square$

Next we characterise those parabolic spheres which arise as Blaschke immersions of special Kähler manifolds.

**Definition 3.4** A parabolic sphere $(M, \nabla, g)$ is called *special* if there exists an almost complex structure $J$ on $M$ such that $g$ is $J$-Hermitian and such that the 2-form $\omega := g(J\cdot, \cdot)$ is $\nabla$-parallel.

**Lemma 3.5** *Under the assumptions of Definition 3.4 the almost complex structure $J$ is integrable and the connection $\nabla^J = J\nabla J^{-1}$ is torsionfree.*

*Proof* Since $\omega$ is Hermitian and $\nabla$-parallel, the same computation as in the proof of Theorem 3.1 shows that $\overline{\nabla} = \nabla^J$. In particular, since, by Theorem 1.1, $\overline{\nabla}$ is torsionfree we conclude that $\nabla^J$ is torsionfree. Of course, since $\nabla$ is flat, $\nabla^J$ is also flat. This shows that the connections $\nabla$ and $\nabla^J$ are both torsionfree and flat. Using this we show that the Nijenhuis tensor

$$N(X,Y) = [JX, JY] - [X,Y] - J[X, JY] - J[JX, Y]$$

for $J$ vanishes. In fact, we may assume that $X$ and $Y$ are $\nabla$-parallel vector fields, and hence $JX$ and $JY$ are $\nabla^J$-parallel. We compute

$$\begin{aligned} N(X,Y) &= -J[X, JY] - J[JX, Y] \\ &= -J(\nabla_X JY - \nabla_{JY} X) - J(\nabla_{JX} Y - \nabla_Y JX) \\ &= \nabla_X^J Y - \nabla_Y^J X = [X,Y] = 0\,. \end{aligned}$$

Now, by the Newlander-Nirenberg Theorem, $J$ is integrable. $\square$

**Theorem 3.6** *A parabolic sphere arises as Blaschke immersion of a special Kähler manifold if and only if it is special.*

*Proof* It is clear that a parabolic sphere which arises as Blaschke immersion of a special Kähler manifold is special in the sense of Definition 3.4. To prove the converse, thanks to Lemma 3.5, we only have to check the condition $d^\nabla J = 0$. But this condition, as already remarked in the proof of Theorem 3.1, is equivalent to the fact that $\nabla^J$ is torsionfree. This was established already in Lemma 3.5. □

We remark that a parabolic 2-sphere is special if and only if its Blaschke metric $g$ is definite. In fact, for a special parabolic sphere, $g$ has the signature of a Hermitian form. Hence, in dimension 2, $g$ is definite. Conversely, any parabolic 2-sphere with definite metric $g$, admits a unique $g$-Hermitian complex structure $J$, as any orientable Riemannian surface does. Then $\omega = g(J\cdot, \cdot)$ coincides with the volume form of $g$, and therefore is $\nabla$-parallel.

The next result shows that the local equations which describe special parabolic spheres are completely integrable. This can be thought of as a generalisation of a classical theorem of Blaschke on parabolic 2-spheres, compare [B, p. 216].

**Corollary 3.7** *Any nondegenerate holomorphic function $F$ of $m$ complex variables defines a special parabolic sphere of dimension $n = 2m$. Conversely, every simply connected special parabolic sphere arises in this way.*